\documentclass[12pt,oneside,reqno]{amsart}
\usepackage{mathrsfs}
\usepackage{graphics}
\usepackage{amssymb}
\usepackage{enumerate}
\newcommand{\dif}{\mathrm{d}}

\newcommand{\be}{\begin{eqnarray}}
\newcommand{\ee}{\end{eqnarray}}
\newcommand{\ce}{\begin{eqnarray*}}
\newcommand{\de}{\end{eqnarray*}}
\newtheorem{theorem}{Theorem}[section]
\newtheorem{lemma}[theorem]{Lemma}
\newtheorem{remark}[theorem]{Remark}
\newtheorem{definition}[theorem]{Definition}
\newtheorem{proposition}[theorem]{Proposition}
\newtheorem{Example}[theorem]{Example}
\newtheorem{corollary}[theorem]{Corollary}
\newtheorem{condition}[theorem]{Condition}
\def\e{\varepsilon}
\def\s{\sigma}
\def\t{\theta}

\def\d{\delta}

\def\g{\gamma}

\def\[{{\Big[}}
\def\]{{\Big]}}
\def\<{{\langle}}
\def\>{{\rangle}}
\def\({{\Big(}}
\def\){{\Big)}}

\def\no{\nonumber}
\def\bt{\begin{theorem}}
\def\et{\end{theorem}}
\def\bl{\begin{lemma}}
\def\el{\end{lemma}}
\def\br{\begin{remark}}
\def\er{\end{remark}}
\def\bx{\begin{Example}}
\def\ex{\end{Example}}
\def\bd{\begin{definition}}
\def\ed{\end{definition}}
\def\bp{\begin{proposition}}
\def\ep{\end{proposition}}
\def\bc{\begin{corollary}}
\def\ec{\end{corollary}}
\def\bco{\begin{condition}}
\def\eco{\end{condition}}

\def\cD{{\mathcal D}}

\def\cK{{\mathcal K}}

\def\cP{{\mathcal P}}

\def\mE{{\mathbb E}}

\def\mN{{\mathbb N}}

\def\mP{{\mathbb P}}

\def\mR{{\mathbb R}}

\def\mW{{\mathbb W}}

\def\sA{{\mathscr A}}
\def\sB{{\mathscr B}}

\def\sL{{\mathscr L}}

\def\sV{{\mathscr V}}

\def\geq{\geqslant}
\def\leq{\leqslant}
\def\epsilon{\varepsilon}
\begin{document}

\allowdisplaybreaks

\title{Limit theorems of invariant measures for multivalued McKean-Vlasov stochastic differential equations}

\author{Huijie Qiao}

\dedicatory{School of Mathematics,
Southeast University,\\
Nanjing, Jiangsu 211189, P.R.China\\
hjqiaogean@seu.edu.cn}

\thanks{{\it AMS Subject Classification(2020):} 60H10}

\thanks{{\it Keywords: Multivalued McKean-Vlasov stochastic differential equations; the exponential ergodicity; the convergence of strong solutions; the convergence of invariant measures}}

\thanks{*This work was partly supported by NSF of China (No.12071071).}

\subjclass{}

\date{}

\begin{abstract}
The work concerns invariant measures for multivalued McKean-Vlasov stochastic differential equations. First of all, we prove the exponential ergodicity of these equations. Then for a sequence of these equations, when their coefficients converge in the suitable sense,  the convergence of corresponding strong solutions are presented. Finally, based on the convergence of these solutions, we establish the convergence of corresponding invariant measures.
\end{abstract}

\maketitle \rm

\section{Introduction}
Given a filtered probability space $(\Omega, \mathscr{F}, \{\mathscr{F}_t\}_{t\in[0,T]}, \mP)$ and a $m$-dimensional Brownian motion $W$ defined on it. Consider the following multivalued McKean-Vlasov stochastic differential equation (SDE for short) on $\mR^d$:
\be
\dif X_t\in \ -A(X_t)\dif t+ \ b(X_t,\sL_{X_t})\dif t+\sigma(X_t,\sL_{X_t})\dif W_t, \quad t\geq 0,
\label{eq1}
\ee
where $A:\mR^d \mapsto 2^{\mR^d} $ is a maximal monotone operator with $\text{Int}(\cD(A))\ne\emptyset$, $\sL_{X_t}:=\mP \circ X_t^{-1}$ is the law of $X_t$ and $X_0$ belongs to $\overline{\cD(A)}$. The coefficients $b:\mR^d\times\cP(\mR^d)\mapsto{\mR^d}, \,\,\sigma:\mR^d\times\cP(\mR^d)\mapsto{\mR^d}\times{\mR^m}$ are Borel measurable, and depend not only on the state $X_{t}$ but also on the distribution $\sL_{X_t}$.

When the coefficients $b,\s$ only depend on the state $X_{t}$, Eq.(\ref{eq1}) is called a multivalued SDE. The study about this equation can be tracked back to Kr\'ee \cite{kp}. From then on, many authors paid attention to these equations (cf. \cite{cepa1, cepa2, cepa3, gz, rwzh,rwzx}). We recall some works. When the operator $A$ is the subdifferential of some convex function, C\'epa and Jacquot \cite{cepa3} proved the ergodicity for the solutions of stochastic variational inequalities by the Bismut formula. Later, Ren, Wu and Zhang \cite{rwzx} studied the exponential ergodicity of multivalued SDEs by proving the irreducibility and strong Feller property. Recently, in \cite{gz}, Guan and Zhang established the weak convergence of invariant measures associated with multivalued SDEs.

When the coefficients $b,\s$ depend on the state $X_{t}$ and the distribution $\sL_{X_t}$, there have been some results about Eq.(\ref{eq1}) (cf. \cite{CHI, G, rw}). When $b,\s$ depend on $\sL_{X_t}$ by integrations, Chi \cite{CHI} showed the well-posedness of Eq.(\ref{eq1}). In the case of $b,\s$ depending on $\sL_{X_t}$ by expectations, Ren and Wang \cite{rw} investigated stochastic variational inequalities, and obtained the well-posedness and a large deviation principle of Eq.(\ref{eq1}). Recently, Gong and the author \cite{G} studied Eq.(\ref{eq1}) under non-Lipschitz conditions. In \cite{G}, they proved  the well-posedness and stability, and presented the existence of invariant measures by Lyapunov functions. Here, we follow up on the result in \cite{G} and give out the exponential ergodicity of Eq.(\ref{eq1}).

Besides, consider the following sequence of multivalued McKean-Vlasov SDEs: for $n\in\mN$,
\be
\dif X^n_t\in \ -A^n(X^n_t)\dif t+b^n(X^n_t,\sL_{X^n_t})\dif t+\sigma^n(X^n_t,\sL_{X^n_t})\dif W_t, \quad t\geq 0,
\label{eqn}
\ee
where $A^n$ is a maximal monotone operator with ${\rm Int}(\cD(A^n))\neq\emptyset$, and $b^n:\mR^d\times\cP(\mR^d)\mapsto{\mR^d}, \,\,\sigma^n:\mR^d\times\cP(\mR^d)\mapsto{\mR^d}\times{\mR^m}$ are Borel measurable. When $b^n\rightarrow b, \sigma^n\rightarrow\sigma$ in some sense as $n\rightarrow\infty$, we will prove that $X_t^n$ also converges to $X_t$ in some suitable sense.

If $A^n=0$ and $b^n, \s^n$ don't depend on the distribution $\sL_{X_t}$, this problem was initially proposed by Stroock and Varadhan in 1979. From then on, similar problems have been discussed by many experts in various formulations (cf. \cite{fi, jjas, ky, q2, q3, sv}). Let us mention some works. Stroock and Varadhan proved the weak convergence in \cite[Chapter 11]{sv}. Later, Kawabata and Yamada \cite{ky} strengthened the conditions and showed that the convergence holds in the $L^1$ sense when $b^n, \s^n$ tended to $b, \s$ respectively in a suitable sense. Recently, Figalli \cite{fi} obtained the weak convergence under bounded conditions. 

If $A^n\neq0$ and $b^n, \s^n$ don't depend on the distribution $\sL_{X_t}$, this problem also has been studied (cf. \cite{ct, gz, rs} and references therein). When the operator $A$ is the subdifferential of the indicator function for a convex set, Rozkosz and Slomi\'nski \cite{rs} proved that the solutions of stochastic variational inequalities weakly converge under bounded conditions. Recently, Guan and Zhang \cite{gz} obtained the mean square convergence for general multivalued SDEs. 

For the case treated here, as far as we know, there seems to be no related result. We establish the mean square convergence for Eq.(\ref{eqn}). Moreover, based on this result, we investigate the convergence of invariant measures for Eq.(\ref{eqn}).

As a whole, our contribution are three-folded:

$\bullet$ We prove the exponential ergodicity of multivalued McKean-Vlasov SDEs.

$\bullet$ For a sequence of multivalued McKean-Vlasov SDEs, we show the stability.

$\bullet$ Based on the above results, we establish the convergence of invariant measures associated with multivalued McKean-Vlasov SDEs.

It is worthwhile to mentioning our methods and results. Since the transition semigroup for Eq.(\ref{eq1}) is not linear, the usual method of proving the ergodicity doesn't work (cf. \cite{q1, rwzx}). Here we show the ergodicity by a contraction estimate. Note that if the operator $A$ is zero, Eq.(\ref{eq1}) is a McKean-Vlasov SDE. Hence, our exponential ergodicity result (Theorem \ref{erg2}) is more general than Theorem 3.1 in \cite{WangFY}. Besides, since multivalued McKean-Vlasov SDEs contain multivalued SDEs, our stability result (Theorem \ref{constrsol}) covers Proposition 2.9 in \cite{gz}.

The paper proceeds as follows. In Section \ref{fram}, we introduce the notions and recall some known results. We place exponential ergodicity of multivalued McKean-Vlasov SDEs in Section \ref{wec}. And in Section \ref{constrsolse} the convergence of strong solutions is presented. Finally, we establish the convergence of invariant measures in Section \ref{coninvmease}.

In the following $C$ with or without indices will denote different positive constants whose values may change from line to line.

\section{Preliminary}\label{fram}

In this section, we will introduce the notions and recall some known results which will be used later in the proofs. 

\subsection{Notations}\label{nota}

In the subsection, we introduce some notations.

Let $\mid\cdot\mid$ and $\|\cdot\|$ be norms of vectors and matrices, respectively. Denote $\langle\cdot$ , $\cdot\rangle$ the scalar product in $\mR^d$ and $B^*$ the transpose of the matrix $B$.

Let $C(\mR^d)$ be the collection of continuous functions on $\mR^d$ and $C_b(\mR^d)$ be the space of all bounded and continuous functions on $\mR^d$. Let $C_{lip,b}(\mR^d)$ be the set of all bounded and Lipschitz continuous functions on $\mR^d$.

Let $\sB(\mR^d)$ be the Borel $\sigma$-algebra on $\mR^d$ and $\cP({\mR^d})$ be the space of all probability measures defined on $\sB(\mR^d)$ carrying the usual topology of the weak convergence. Define  the following metric on $\cP(\mR^d)$:
\ce
\rho(\mu,\nu):=\sup\limits_{\parallel{\varphi}\parallel_{C_{lip,b}(\mR^d)\leq1}}\left|{\int_{\mR^d}\varphi(x)\mu(\dif x)-\int_{\mR^d}\varphi(x)\nu(\dif x)}\right|, \quad \mu, \nu\in\cP(\mR^d),
 \de
where 
$$
\parallel{\varphi}\parallel_{C_{lip,b}(\mR^d)}:=\sup\limits_{x\in\mR^d}|\varphi(x)|+\sup\limits_{x\neq y}\frac{|\varphi(x)-\varphi(y)|}{|x-y|}.
$$
And $(\cP({\mR^d}),\rho)$ is a Polish space. Moreover, the topology generated by $\rho$ on $\cP(\mR^d)$ coincides with the usual topology of weak convergence (\cite[Theorem 3.1]{kr}.

For any $\t\in[1,\infty)$, we consider the following subspace of $\cP(\mR^d)$:
$$
\cP_\t\left( \mathbb{R}^d \right) :=\left\{ \mu \in \cP\left( \mathbb{R}^d \right): \|\mu\|_\t^\t:=\int_{\mathbb{R}^d}{\left| x \right|^\t\mu \left( \dif x \right) <\infty} \right\}.
$$
It is known that $\cP_\t(\mR^d)$ is a Polish space endowed with the Wasserstein distance defined by
$$
\mathbb{W}_\t(\mu,\nu):= \inf\limits_{\pi\in\Psi(\mu,\nu)}\left(\int_{\mathbb{R}^d\times\mathbb{R}^d}|x-y|^\t\pi(\dif x,\dif y)\right)^{\frac{1}{\t}}, \quad \mu , \nu\in \cP_\t(\mR^d),
$$
where $\Psi(\mu,\nu)$ is the set of all couplings $\pi$ with marginal distributions $\mu$ and $\nu$. Furthermore, when $\t=1$, $\rho$ is equivalent to $\mW_1$ on $\cP_1(\mR^d)$  (\cite[Theorem 8.10.45, P. 235]{bo}).

\subsection{Maximal monotone operators}\label{mmo}

In the subsection, we introduce maximal monotone operators.

For a multivalued operator $A: \mR^d\mapsto 2^{\mR^d}$, where $2^{\mR^d}$ stands for all the subsets of $\mR^d$, set
\ce
&&\cD(A):= \left\{x\in \mR^d: A(x) \ne \emptyset\right\},\\
&&Gr(A):= \left\{(x,y)\in \mR^{2d}:x \in \cD(A), ~ y\in A(x)\right\}.
\de
We say that $A$ is monotone if $\langle x_1 - x_2, y_1 - y_2 \rangle \geq 0$ for any $(x_1,y_1), (x_2,y_2) \in Gr(A)$, and $A$ is maximal monotone if 
$$
(x_1,y_1) \in Gr(A) \iff \langle x_1-x_2, y_1 -y_2 \rangle \geq 0, \quad \forall (x_2,y_2) \in Gr(A).
$$

In the following, we recall some properties of a maximal monotone operator $A$ (cf.\cite{cepa1, gp}):
\begin{enumerate}[(i)]
\item
${\rm Int}(\cD(A))$ and $\overline{\mathrm{\cD}(A)}$ are convex subsets of $\mR^d$ with ${\rm Int}\left( \overline{\mathrm{\cD}(A)} \right) = {\rm Int}\( \mathrm{\cD}(A) \) 
$, where ${\rm Int}(\cD(A))$ denotes the interior of the set $\cD(A)$. 
\item For every $x\in\mR^d$, $A(x)$ is a closed and convex subset of $\mR^d$. Let $A^{\circ}\left( x \right):= proj_{A(x)}(0)$ be the minimal section of $A$, where $proj_D$ is designated as the projection on every closed and
convex subset $D\subset\mR^d$ and $proj_{\emptyset}(0) =\infty$. Then
$$
x\in\cD(A) \Longleftrightarrow |A^{\circ}\left( x \right)|<\infty.
$$
\item For $\varepsilon > 0$, the resolvent operator $J_{\varepsilon}:=\left( I+\varepsilon A \right) ^{-1}$ is a single-valued and contractive operator defined on $\mR^d$ and takes values in $\mathrm{\cD}(A)$, and 
$$
\lim\limits_{\e\downarrow 0}J_\e(x)=proj_{\overline{\mathrm{\cD}(A)}}(x), \quad x\in\mR^d.
$$
\item $A_{\varepsilon}:= \frac{1}{\varepsilon}\left(I-J_{\varepsilon}\right)$, called the Yosida approximation of $A$, is also a single-valued, maximal monotone and Lipschitz continuous operator with the Lipschitz constant $\frac{1}{\varepsilon}$. 
\item 
$A_{\varepsilon}\left( x \right) \,\,\in \,\,A\left( J_{\varepsilon}\left( x \right) \right), \quad x\in\mR^d$.
\item 
$|A_{\varepsilon}\left( x \right) | \leq |A^{\circ}\left( x \right)|, \quad x\in\cD(A)$.
\item
For any $x\in\cD(A)$, $\lim\limits_{\varepsilon \downarrow 0}A_{\varepsilon}\left( x \right) =A^{\circ}\left( x \right)$, and
\ce
&&\lim_{\varepsilon \downarrow 0}|A_{\varepsilon}\left( x \right) |=|A^{\circ}\left( x \right)|, \quad x\in\cD(A),\\
&&\lim_{\varepsilon \downarrow 0}|A_{\varepsilon}\left( x \right) |=\infty,\qquad\qquad x\notin\cD(A).
\de
\end{enumerate}

About the Yosida approximation $A_{\varepsilon}$, we also mention the following property (\cite[Lemma 5.4]{cepa1}). 

\bl\label{yosi}
There exist three constants $a\in\mR^d, M_1>0, M_2\geq 0$ only dependent on $A$ such that for any $\e>0$ and $x\in\mR^d$
$$
\<A_{\varepsilon}(x), x-a\>\geq M_1|A_{\varepsilon}(x)|-M_2|x-a|-M_1M_2.
$$
\el

\medspace

Take any $T>0$ and fix it. Let $\sV_{0}$ be the set of all continuous functions $K: [0,T]\mapsto\mR^d$ with finite variations and $K_{0} = 0$. For $K\in\sV_0$ and $s\in [0,T]$, we shall use $|K|_{0}^{s}$ to denote the variation of $K$ on [0,s]
and write $|K|_{TV}:=|K|_{0}^{T}$. Set
\ce
&&\sA:=\Big\{(X,K): X\in C([0,T],\overline{\cD(A)}), K \in \sV_0, \\
&&\qquad\qquad\quad~\mbox{and}~\langle X_{t}-x, \dif K_{t}-y\dif t\rangle \geq 0 ~\mbox{for any}~ (x,y)\in Gr(A)\Big\}.
\de
And about $\sA$ we have the following result (cf.\cite{cepa2, ZXCH}).

\bl\label{equi}
For $X\in C([0,T],\overline{\cD(A)})$ and $K\in \sV_{0}$, the following statements are equivalent:
\begin{enumerate}[(i)]
	\item $(X,K)\in \sA$.
	\item For any $(x,y)\in C([0,T],\mR^d)$ with $(x_t, y_t)\in Gr(A)$, it holds that 
	$$
	\left\langle X_t-x_t, \dif K_t-y_t\dif t\right\rangle \geq0.
	$$
	\item For any $(X^{'},K^{'})\in \sA$, it holds that 
	$$
	\left\langle X_t-X_t^{'},\dif K_t-\dif K_t^{'}\right\rangle \geq0.
	$$
\end{enumerate}
\el

\subsection{Multivalued McKean-Vlasov SDEs}

In the subsection, we introduce multivalued McKean-Vlasov SDEs.

First of all, we define strong solutions for Eq.$(\ref{eq1})$. 

\bd\label{strosolu}
We say that Eq.$(\ref{eq1})$ admits a strong solution with the initial value $X_0$ if there exists a pair of adapted processes $(X,K)$ on a filtered probability space $(\Omega, \mathscr{F}, \{\mathscr{F}_t\}_{t\in[0,T]}, \mP)$ such that

(i) $X_t\in{\mathscr{F}_t^W}$, where $\{\mathscr{F}_t^W\}_{t\in[0,T]}$ stands for the $\sigma$-field filtration generated by $W$,

(ii) $(X_{\cdot}(\omega),K_{\cdot}(\omega))\in \sA$ a.s. $\mP$,

(iii) it holds that
\ce
\mP\left\{\int_0^T(\mid{b(X_s,\sL_{X_s})}\mid+\parallel{\sigma(X_s,\sL_{X_s})}\parallel^2)\dif s<+\infty\right\}=1,
\de
and
\ce
X_t=X_0-K_{t}+\int_0^tb(X_s,\sL_{X_s})\dif s+\int_0^t\sigma(X_s,\sL_{X_s})\dif W_s, \quad 0\leq{t}\leq{T}.
\de
\ed

Next, if Eq.(\ref{eq1}) has a unique strong solution $(X_{\cdot},K_{\cdot})$, we define $P^*_{t}\mu:=\sL_{X_t}$ for $\sL_{X_0}=\mu\in\cP_2(\overline{\cD(A)})$.

\bd\label{iae}
We call $\mu\in\cP_2(\overline{\cD(A)})$ an invariant probability measure of $\{P_{t}^{*}, t\geq 0\}$ if $P_{t}^{*}\mu=\mu$ for all $t\geq0$. We call $\{P_{t}^{*}, t\geq 0\}$ ergodic if there exists $\mu\in\cP_2(\overline{\cD(A)})$ such that $\lim\limits_{t\rightarrow\infty}P_{t}^{*}\nu=\mu$ weakly for any $\nu\in\cP_2(\overline{\cD(A)})$.
\ed
It is obvious that an ergodic $\{P_{t}^{*}, t\geq 0\}$ has a unique invariant probability measure.

\section{The exponential ergodicity}\label{wec}

In this section we are devoted to studying the exponential ergodicity for multivalued McKean-Vlasov SDEs. 

We make the following assumptions:
\begin{enumerate}[(${\bf H}_{A}$)]
\item $0\in \cD(A)$.
\end{enumerate}
\begin{enumerate}[(${\bf H}^1_{b,\s}$)]
	\item The function $b$ is continuous in $(x,\mu)$, and there exists a constant $L_{b,\s}>0$ such that for any $(x,\mu)\in\mR^{d}\times{\cP_{2}(\mR^d)}$
	\ce
	|{b(x,\mu)}|+\|\s(x,\mu)\|\leq L_{b,\s}(1+|x|+\|\mu\|_2).
	\de
\end{enumerate}
\begin{enumerate}[(${\bf H}^2_{b}$)]
	\item There exists a constant $L_1>0$ such that for any $(x_{1},\mu_{1}),(x_{2},\mu_{2})\in\mR^{d}\times\cP_{2}(\mR^d)$
	\ce
	2\langle x_{1}-x_{2}, b(x_{1},\mu_{1})-b(x_{2},\mu_{2})\rangle\leq L_1(|x_{1}-x_{2}|^{2}+\mW^2_{2}(\mu_1,\mu_2)),
	\de
\end{enumerate}
\begin{enumerate}[(${\bf H}^2_{\s}$)]
	\item There exists a constant $L_2>0$ such that for any $(x_{1},\mu_{1}),(x_{2},\mu_{2})\in\mR^{d}\times\cP_{2}(\mR^d)$
	\ce
	\|\sigma(x_{1},\mu_{1})-\sigma(x_{2},\mu_{2})\|^2\leq L_2(|x_{1}-x_{2}|^{2}+\mW^2_{2}(\mu_1,\mu_2)).
	\de
\end{enumerate}
\begin{enumerate}[(${\bf H}^{2'}_{b,\s}$)]
	\item There exist two constants $L_3, L_4>0, L_4-L_3>2L_2$ such that for any  $(x_{1},\mu_{1}),(x_{2},\mu_{2})\in\mR^{d}\times\cP_{2}(\mR^d)$
	\ce
	2\langle x_{1}-x_{2}, b(x_{1},\mu_{1})-b(x_{2},\mu_{2})\rangle+\|\sigma(x_{1},\mu_{1})-\sigma(x_{2},\mu_{2})\|^2\leq L_3\mW^{2}_2(\mu_{1},\mu_{2})-L_4|x_{1}-x_{2}|^{2}.
	\de
\end{enumerate}

\br
(${\bf H}^1_{b,\s}$) (${\bf H}^2_{b}$) (${\bf H}^2_{\s}$) assure the well-posedness of Eq.(\ref{eq1}) (cf. \cite[Theorem 3.5]{G}). Moreover, (${\bf H}^{2'}_{b,\s}$) implies (${\bf H}^2_{b}$).
\er

Under the above conditions, we have the following conclusion which is the main result in this section.

\bt\label{erg2}
Assume that (${\bf H}_{A}$) (${\bf H}^1_{b,\s}$) (${\bf H}^2_{\s}$) (${\bf H}^{2'}_{b,\s}$) hold and $\mE|X_0|^2<\infty$. Then $\{P_{t}^{*}, t\geq 0\}$ has a unique invariant probability measure $\mu_{\infty}\in\cP_{2}(\overline{\cD(A)})$ such that for any $\nu_{0}\in\cP_{2}(\overline{\cD(A)})$
\be
\mW^{2}_{2}(P_{t}^{*}\nu_{0},\mu_{\infty})\leq 2(\|\nu_0\|^2_2+\|\mu_{\infty}\|^2_2)e^{-\lambda t}, \quad t\geq 0,
\label{distesti}
\ee
where $\lambda:=L_4-L_3$.
\et

To prove the above theorem, we prepare a key proposition.

\bp\label{erg1}
Assume that (${\bf H}^1_{b,\s}$) (${\bf H}^2_{\s}$) (${\bf H}^{2'}_{b,\s}$) hold. For any $\mu_{0},\nu_{0}\in\cP_{2}(\overline{\cD(A)})$,
\be
\mW^2_{2}(\sL_{X_t},\sL_{Y_t})\leq \mW^{2}_{2}(\mu_{0},\nu_{0})e^{-\lambda t}, \quad t\geq 0,
\label{xtytdiff}
\ee
where $(X,K), (Y,\tilde{K})$ are two solutions to Eq.(\ref{eq1}) such that $\sL_{X_{0}}=\mu_{0},\sL_{Y_{0}}=\nu_{0}$ and
$$
\mW^{2}_{2}(\mu_{0},\nu_{0})=\mE|X_{0}-Y_{0}|^{2}.
$$
\ep
\begin{proof}
First of all, since $(X,K), (Y,\tilde{K})$ are two solutions to Eq.(\ref{eq1}), it holds that
\ce
&&X_{t}=X_{0}-K_{t}+\int_{0}^{t}b(X_{s},\sL_{X_s})\dif s+\int_{0}^{t}\sigma(X_{s},\sL_{X_s})\dif W_{s},\\
&&Y_{t}=Y_{0}-\tilde{K}_{t}+\int_{0}^{t}b(Y_{s},\sL_{Y_s})\dif s+\int_{0}^{t}\sigma(Y_{s},\sL_{Y_s})\dif W_{s}.
\de
By  It\^o's formula and Lemma \ref{equi}, we have
\ce
|X_{t}-Y_{t}|^{2}e^{\lambda t}&=&|X_{0}-Y_{0}|^{2}-2\int_{0}^{t}e^{\lambda s}\left\langle X_{s}-Y_{s},\dif K_{s}-\dif\tilde{K}_{s}  \right\rangle\\
&&+2\int_{0}^{t}e^{\lambda s}\left\langle X_{s}-Y_{s},b(X_{s},\sL_{X_s})-b(Y_{s},\sL_{Y_s}) \right\rangle \dif s\\
&&+2\int_{0}^{t}e^{\lambda s}\left\langle X_{s}-Y_{s},[\sigma(X_{s},\sL_{X_s})-\sigma(Y_{s},\sL_{Y_s}) ]\dif W_{s}\right\rangle\\
&&+\int_{0}^{t}e^{\lambda s} \|\sigma(X_{s},\sL_{X_s})-\sigma(Y_{s},\sL_{Y_s})\|^{2} \dif s+\int_0^t\lambda e^{\lambda s}|X_{s}-Y_{s}|^{2}\dif s\\
&\leq&|X_{0}-Y_{0}|^{2}+2\int_{0}^{t}e^{\lambda s}\left\langle X_{s}-Y_{s},b(X_{s},\sL_{X_s})-b(Y_{s},\sL_{Y_s}) \right\rangle \dif s\\
&&+2\int_{0}^{t}e^{\lambda s}\left\langle X_{s}-Y_{s},[\sigma(X_{s},\sL_{X_s})-\sigma(Y_{s},\sL_{Y_s}) ]\dif W_{s}\right\rangle\\
&&+\int_{0}^{t} e^{\lambda s}\|\sigma(X_{s},\sL_{X_s})-\sigma(Y_{s},\sL_{Y_s})\|^{2} \dif s+\int_0^t\lambda e^{\lambda s}|X_{s}-Y_{s}|^{2}\dif s.
\de
Taking the expectation on two sides, by (${\bf H}^{2'}_{b,\s}$) one can obtain that
\be
\mE|X_{t}-Y_{t}|^{2}e^{\lambda t}&\leq&\mE|X_{0}-Y_{0}|^{2}+\int_{0}^{t}e^{\lambda s}\[L_3\mW_2^2(\sL_{X_s}, \sL_{Y_s})-L_4\mE|X_{s}-Y_{s}|^{2}\]\dif s\no\\
&&+\int_0^t\lambda e^{\lambda s}\mE|X_{s}-Y_{s}|^{2}\dif s\no\\
&\leq&\mE|X_{0}-Y_{0}|^{2}-\lambda \int_{0}^{t}e^{\lambda s}\mE|X_{s}-Y_{s}|^{2}\dif s+\int_0^t\lambda e^{\lambda s}\mE|X_{s}-Y_{s}|^{2}\dif s\no\\
&=&\mW^{2}_{2}(\mu_{0},\nu_{0}),
\label{xynumu}
\ee
where the fact that $\mW^2_{2}(\sL_{X_s},\sL_{Y_s})\leq\mE|X_{s}-Y_{s}|^{2}$ is used. Therefore, (\ref{xtytdiff}) holds.
\end{proof}

Now, we prove Theorem \ref{erg2}.

{\bf Proof of Theorem \ref{erg2}.}
{\bf Step 1.} Let $(X^0,K^0)$ be the solution to Eq.(\ref{eq1}) with $\sL_{X_0^0}=\d_0$, where $\d_0$ is the Dirac measure in the point $0$. Then we prove that there exists a probability measure $\mu_{\infty}\in\cP_2(\mR^d)$ such that 
\be
\lim\limits_{t\rightarrow\infty}\mW_2(\sL_{X_t^0},\mu_{\infty})=0.
\label{findlimi}
\ee

First of all, by the uniqueness in law of weak solutions for Eq.(\ref{eq1}), it holds that
\be
\sL_{X_{t+s}^0}=\sL_{X_{t}^{X_s^0}}.
\label{uniqweak}
\ee
Thus, (\ref{xtytdiff}) implies that
\ce
\lim\limits_{t\rightarrow\infty}\sup\limits_{s\geq 0}\mW_2^2(\sL_{X_t^0},\sL_{X_{t+s}^0})=\lim\limits_{t\rightarrow\infty}\sup\limits_{s\geq 0}\mW_2^2(\sL_{X_t^0},\sL_{X_{t}^{X_s^0}})\leq \lim\limits_{t\rightarrow\infty}\sup\limits_{s\geq 0}\mE|X_s^0|^2e^{-\lambda t}.
\de
We claim that 
\be
\sup\limits_{s\geq 0}\mE|X_s^0|^2<\infty.
\label{xs0}
\ee
So, $\{\sL_{X_t^0}: t\geq 0\}$ is a Cauchy sequence in $(\cP_2(\mR^d), \mW_2)$ and there exists a probability measure $\mu_{\infty}\in\cP_2(\mR^d)$ satisfying (\ref{findlimi}).

Here we justify the claim (\ref{xs0}). Note that (${\bf H}^{2}_{\sigma}$) (${\bf H}^{2'}_{b,\sigma}$) give the following inequality
\be
2\<x, b(x,\mu)\>+\|\sigma(x,\mu)\|^2\leq C+\(L_3+L_2\)\|\mu\|_2^2-\(L_4-L_2-\frac{\lambda-2L_2}{4}\)|x|^2.
\label{ergcon}
\ee
Thus, fixing $y\in A(0)$ and $0<\eta<\frac{\lambda-2L_2}{2}$, applying the It\^o formula to $|X_s^0|^2e^{\eta s}$ and taking the expectation on two sides, by Lemma \ref{equi} and the Young inequality we get that
\ce
\mE|X_s^0|^2e^{\eta s}&=&\mE\int_0^s\eta e^{\eta r}|X_r^0|^2\dif r-2\mE\int_0^se^{\eta r}\<X_r^0, \dif K_r^0\>\\
&&+\mE\int_0^se^{\eta r}\left[2\<X_r^0, b(X_r^0, \sL_{X_r^0})\>+\|\sigma(X_r^0, \sL_{X_r^0})\|^2\right]\dif r\\
&\leq&\mE\int_0^s\eta e^{\eta r}|X_r^0|^2\dif r+2\mE\int_0^se^{\eta r}|y||X_r^0|\dif r+C\int_0^se^{\eta r}\dif r\\
&&+\mE\int_0^se^{\eta r}\left[\(L_3+L_2\)\|\sL_{X_r^0}\|_2^2-\(L_4-L_2-\frac{\lambda-2L_2}{4}\)|X_r^0|^2\right]\dif r\\
&\leq&C\int_0^se^{\eta r}\dif r+\left(\eta+\frac{\lambda-2L_2}{4}\right)\mE\int_0^s e^{\eta r}|X_r^0|^2\dif r-\frac{3(\lambda-2L_2)}{4}\mE\int_0^s e^{\eta r}|X_r^0|^2\dif r\\
&\leq&C\frac{e^{\eta s}-1}{\eta},
\de
where the fact $\|\sL_{X_r^0}\|_2^2=\mE|X_r^0|^2$ is used. That is, (\ref{xs0}) is right.

{\bf Step 2.} We prove (\ref{distesti}).

Note that for any $s\geq 0$
\ce
\mW_2(\sL_{X_s^{\mu_{\infty}}}, \mu_{\infty})&\leq& \mW_2(\sL_{X_s^{\mu_{\infty}}}, \sL_{X_s^{X_t^0}})+\mW_2(\sL_{X_s^{X_t^0}}, \sL_{X_t^0})+\mW_2(\sL_{X_t^0}, \mu_{\infty})\\
&\overset{(\ref{xtytdiff})}{\leq}&\mW_2(\mu_{\infty}, \sL_{X_t^0})e^{-\lambda s/2}+\mW_2(\sL_{X_s^0}, \d_0)e^{-\lambda t/2}+\mW_2(\sL_{X_t^0}, \mu_{\infty}).
\de
So, this together with (\ref{findlimi}) yields that $\sL_{X_s^{\mu_{\infty}}}=\mu_{\infty}$ and 
\ce
\mW^2_2(\sL_{X_t^{\nu_0}}, \mu_{\infty})=\mW^2_2(\sL_{X_t^{\nu_0}}, \sL_{X_t^{\mu_{\infty}}})\overset{(\ref{xtytdiff})}{\leq}\mE|X_0^{\nu_0}-X_0^{\mu_{\infty}}|^2e^{-\lambda t}\leq 2(\|\nu_0\|^2_2+\|\mu_{\infty}\|^2_2)e^{-\lambda t},
\de
which completes the proof.

\section{Convergence of strong solutions}\label{constrsolse}

In the section, we observe the convergence of strong solutions for a sequence of multivalued McKean-Vlasov SDEs.

Consider Eq.(\ref{eq1}) and Eq.(\ref{eqn}), i.e.
\ce
&&\dif X_t\in \ -A(X_t)\dif t+ \ b(X_t,\sL_{X_t})\dif t+\sigma(X_t,\sL_{X_t})\dif W_t, \qquad \quad t\geq 0,\\
&&\dif X^n_t\in \ -A^n(X^n_t)\dif t+b^n(X^n_t,\sL_{X^n_t})\dif t+\sigma^n(X^n_t,\sL_{X^n_t})\dif W_t, \quad t\geq 0.
\de
When $b, \s, b^n, \s^n$ uniformly satisfy (${\bf H}^1_{b,\s}$) (${\bf H}^2_{b}$) (${\bf H}^2_{\s}$), by \cite[Theorem 3.5]{G}, Eq.(\ref{eq1}) with $X_0\in\overline{\cD(A)}, \mE|X_0|^2<\infty$ and Eq.(\ref{eqn}) with $X^n_0\in\overline{\cD(A^n)}, \mE|X^n_0|^2<\infty$ have unique solutions $(X_{\cdot}, K_{\cdot}), (X^n_{\cdot}, K^n_{\cdot})$, respectively.

We also assume:
\begin{enumerate}[(${\bf H}^1_{A^n,A}$)]
\item $\cD(A)=\cD(A^n)$, $0\in {\rm Int}(\cD(A))$ and $A^n$ is locally bounded at $0$ uniformly in $n$, i.e., there exists $\kappa>0$ such that 
\ce
\g:=\sup\limits_{n}\sup\{|y|; y\in A^n(x), x\in B(0,\kappa):=\{x\in\mR^d: |x|\leq \kappa\}\subset\cD(A^n)\}<\infty.
\de
\end{enumerate}
\begin{enumerate}[(${\bf H}^2_{A^n,A}$)]
\item For any $\e>0$ and any compact set $K\subset \overline{\cD(A)}$, 
$$
\lim\limits_{n\rightarrow\infty}\sup\limits_{x\in K}|A_{\e}^n x-A_{\e}x|=0,
$$
where $A_{\e}^n, A_{\e}$ are the Yosida approximation of $A^n, A$, respectively.
\end{enumerate}
\begin{enumerate}[(${\bf H}_{b^n,b,\s^n,\s}$)]
\item For any compact set $\cK\in\overline{\cD(A)}\times\cP_2(\overline{\cD(A)})$, 
\ce
\lim\limits_{n\rightarrow\infty}\sup\limits_{(x,\mu)\in\cK}|b^n(x,\mu)-b(x,\mu)|=0, \quad \lim\limits_{n\rightarrow\infty}\sup\limits_{(x,\mu)\in\cK}\|\s^n(x,\mu)-\s(x,\mu)\|=0.
\de
 \end{enumerate}
 
Here we characterize the relationship between $X^n$ and $X$ in the following theorem which is the main result in this section.

\bt\label{constrsol}
Assume that $b, \s, b^n, \s^n$ uniformly satisfy (${\bf H}^1_{b,\s}$) (${\bf H}^2_{b}$) (${\bf H}^2_{\s}$) and $X_0=X^n_0, \mE|X_0|^{2p}<\infty$ for any $p>2$. Moreover, (${\bf H}^1_{A^n,A}$) (${\bf H}^2_{A^n,A}$) (${\bf H}_{b^n,b,\s^n,\s}$) hold. Then for any $T>0$
\ce
\lim\limits_{n\rightarrow\infty}\sup\limits_{t\in[0,T]}\mE|X^n_t-X_t|^2=0.
\de
\et

To prove the above theorem, we prepare some key lemmas and propositions. The following lemma is a more general version of \cite[Corollary 3.6]{G}.

\bl\label{xkes}
Suppose that $b,\s, b^n, \s^n$ uniformly satisfy (${\bf H}^1_{b,\s}$) (${\bf H}^2_{b}$) (${\bf H}^2_{\s}$) and $\mE|X_0|^{2p}<\infty, \sup\limits_{n}\mE|X^n_0|^{2p}<\infty$ for any $p\geq 1$. Then we have that 
\ce
\mE\left(\sup\limits_{t\in[0,T]}|X_t|^{2p}\right)+\mE|K|^T_0\leq C, \quad \sup\limits_{n}\mE\left(\sup\limits_{t\in[0,T]}|X^n_t|^{2p}\right)+\sup\limits_{n}\mE|K^n|^T_0\leq C.
\de
\el

Since the proof of the above lemma is similar to that for \cite[Corollary 3.6]{G}, we omit it.

Next, for any $\varepsilon>0$, consider the penalized versions of Eq.(\ref{eq1}) and Eq.(\ref{eqn}) :
\be
&&\dif X^\e_t=\ -A_{\e}(X^\e_t)\dif t+ \ b(X^\e_t,\sL_{X^\e_t})\dif t+\sigma(X^\e_t,\sL_{X^\e_t})\dif W_t,\label{Yeq1}\\
&&\dif X^{n,\e}_t=\ -A^n_{\e}(X^{n,\e}_t)\dif t+ \ b^n(X^{n,\e}_t,\sL_{X^{n,\e}_t})\dif t+\sigma^n(X^{n,\e}_t,\sL_{X^{n,\e}_t})\dif W_t.
\label{Yeq2}
\ee
Note that $A_{\varepsilon}, A^n_{\e}$ are single-valued, maximal monotone and Lipschitz continuous functions (cf. Subsection \ref{mmo}). Thus, by \cite[Theorem 3.1]{DQ1} or \cite[Theorem 2.1]{WangFY}, we know that under (${\bf H}^1_{b,\s}$) (${\bf H}^2_{b}$) (${\bf H}^2_{\s}$) Eq.(\ref{Yeq1}) with $X^\e_0=X_0$ and Eq.(\ref{Yeq2}) with $X^{n,\e}_0=X^n_0$ have unique strong solutions $X^\e, X^{n,\e}$, respectively. Moreover, we give some uniform moment estimates about $X^\e, X^{n,\e}$. 

\bl\label{Yeqsomo}
Suppose that $b,\s, b^n, \s^n$ uniformly satisfy (${\bf H}^1_{b,\s}$) (${\bf H}^2_{b}$) (${\bf H}^2_{\s}$) and $\mE|X_0|^{2p}<\infty, \sup\limits_{n}\mE|X^n_0|^{2p}<\infty$ for any $p\geq 1$. Assume that (${\bf H}^1_{A^n,A}$) hold. Then we have that 
\ce
\sup\limits_{\e}\mE\left(\sup\limits_{t\in[0,T]}|X_t^\e|^{2p}\right)\leq C, \quad \sup\limits_{n,\e}\mE\left(\sup\limits_{t\in[0,T]}|X_t^{n,\e}|^{2p}\right)\leq C.
\de
\el
\begin{proof}
First of all, by the It\^o formula, it holds that for $t\in[0,T]$
\ce
|X_t^\e-a|^{2p}&=&|X_0-a|^{2p}-\int_0^t 2p|X_s^\e-a|^{2p-2}\<X_s^\e-a,A_{\e}(X^\e_s)\>\dif s\no\\
&&+\int_0^t 2p|X_s^\e-a|^{2p-2}\<X_s^\e-a,b(X^\e_s,\sL_{X^\e_s})\>\dif s\no\\
&&+\int_0^t 2p|X_s^\e-a|^{2p-2}\<X_s^\e-a,\sigma(X^\e_s,\sL_{X^\e_s})\dif W_s\>\no\\
&&+\int_0^t2p(p-1)|X_s^\e-a|^{2p-4}(X_s^\e-a)^*\sigma\s^*(X^\e_s,\sL_{X^\e_s})(X_s^\e-a)\dif s\no\\
&&+\int_0^tp|X_s^\e-a|^{2p-2}\|\sigma(X^\e_s,\sL_{X^\e_s})\|^2\dif s,
\de
where $a\in\mR^d$ is the same to that in Lemma \ref{yosi}. And Lemma \ref{yosi} and (${\bf H}^1_{b,\s}$) imply that
\ce
|X_t^\e-a|^{2p}&\leq&|X_0-a|^{2p}-2pM_1\int_0^t|X_s^\e-a|^{2p-2}|A_{\e}(X^\e_s)|\dif s+2pM_2\int_0^t|X_s^\e-a|^{2p-1}\dif s\no\\
&&+2M_1M_2p\int_0^t|X_s^\e-a|^{2p-2}\dif s+\int_0^tp|X_s^\e-a|^{2p}\dif s\no\\
&&+\int_0^tp|X_s^\e-a|^{2p-2}|b(X^\e_s,\sL_{X^\e_s})|^2\dif s+\int_0^t2p^2|X_s^\e-a|^{2p-2}\|\sigma(X^\e_s,\sL_{X^\e_s})\|^2\dif s\no\\
&&+\int_0^t 2p|X_s^\e-a|^{2p-2}\<X_s^\e-a,\sigma(X^\e_s,\sL_{X^\e_s})\dif W_s\>\no\\
&\leq&|X_0-a|^{2p}+(2M_2+2M_1M_2+1)p\int_0^t|X_s^\e-a|^{2p}\dif s+2(M_1+1)M_2pt\no\\
&&+3(2p^2+p)L^2_{b,\s}\int_0^t|X_s^\e-a|^{2p-2}(1+|X^\e_s|^2+\|\sL_{X^\e_s}\|_2^2)\dif s\no\\
&&+\int_0^t 2p|X_s^\e-a|^{2p-2}\<X_s^\e-a,\sigma(X^\e_s,\sL_{X^\e_s})\dif W_s\>.
\de
By the BDG inequality and the Young inequality, we get that
\ce
&&\mE\left(\sup\limits_{t\in[0,T]}|X_t^\e-a|^{2p}\right)\\
&\leq&\mE|X_0-a|^{2p}+(2M_2+2M_1M_2+1)p\int_0^T\mE|X_s^\e-a|^{2p}\dif s+2(M_1+1)M_2pT\no\\
&&+C\int_0^T\mE|X_s^\e-a|^{2p}\dif s+C\int_0^T\mE(1+|X^\e_s|^{2p}+\mE|X^\e_s|^{2p})\dif s\no\\
&&+C\mE\left(\int_0^T |X_s^\e-a|^{4p-2}\|\sigma(X^\e_s,\sL_{X^\e_s})\|^2\dif s\right)^{1/2}\no\\
&\leq&\mE|X_0-a|^{2p}+(2M_2+2M_1M_2+1)p\int_0^T\mE|X_s^\e-a|^{2p}\dif s+2(M_1+1)M_2pT\no\\
&&+C\int_0^T\mE|X_s^\e-a|^{2p}\dif s+C\int_0^T\mE(1+|X^\e_s|^{2p}+\mE|X^\e_s|^{2p})\dif s\no\\
&&+\frac{1}{2}\mE\left(\sup\limits_{t\in[0,T]}|X_t^\e-a|^{2p}\right)\\
&\leq&C+C\int_0^T\mE\left(\sup\limits_{r\in[0,s]}|X_r^\e-a|^{2p}\right)\dif s+\frac{1}{2}\mE\left(\sup\limits_{t\in[0,T]}|X_t^\e-a|^{2p}\right),
\de
which together with the Gronwall inequality yields that
\ce
\mE\left(\sup\limits_{t\in[0,T]}|X_t^\e-a|^{2p}\right)\leq Ce^{CT}.
\de

Finally, note that for any $x\in\cD(A)$
\be
-\<x,A^n_\e(x)\>=-\<x,A^n_\e(x)-A^n_\e(0)+A^n_\e(0)\>\leq -\<x,A^n_\e(0)\>\leq \frac{1}{2}|x|^2+\frac{1}{2}\g^2,
\label{anees}
\ee
where the monotonicity of $A^n_\e$ and the fact that $\sup\limits_{n}|A^n_\e(0)|\leq\sup\limits_{n}|(A^n)^\circ(0)|\leq \g$ are used. So, taking $a=0$ and replacing Lemma \ref{yosi} by (\ref{anees}), by the same deduction to that of the first estimate, we obtain the second estimate. The proof is complete.
\end{proof}

\bp\label{xexcon}
Suppose that $b,\s, b^n, \s^n$ uniformly satisfy (${\bf H}^1_{b,\s}$) (${\bf H}^2_{b}$) (${\bf H}^2_{\s}$) and $X_0=X^n_0, \mE|X_0|^{2p}<\infty$ for any $p>2$. Assume that (${\bf H}^1_{A^n,A}$) hold. Then we have that
\ce
\lim\limits_{\e\rightarrow0}\mE\left(\sup\limits_{t\in[0,T]}|X_t^\e-X_t|^2\right)=0,\quad \lim\limits_{\e\rightarrow0}\sup\limits_{n}\mE\left(\sup\limits_{t\in[0,T]}|X_t^{n,\e}-X^n_t|^2\right)=0.
\de
\ep
\begin{proof}
First of all, for any $t\in[0,T]$, set  
$$
M_t:=\int_0^t b(X_s,\sL_{X_s})\dif s+\int_0^t\sigma(X_s,\sL_{X_s})\dif W_s, 
$$
and then by the H\"older inequality, the BDG inequality, (${\bf H}^1_{b,\s}$) and Lemma \ref{xkes} it holds that
\be
\mE\left(\sup\limits_{t\in[0,T]}|M_t|^{2p}\right)\leq C, \quad \mE\left(\sup\limits_{|t-s|\leq \d}|M_t-M_s|^{2p}\right)\leq C\d^{p-1}.
\label{mcont}
\ee
So, we construct the following equation on $\mR^d$:
\ce\left\{\begin{array}{l}
\dif \bar{X}_t^\e=-A_\e(\bar{X}_t^\e)\dif t+\dif M_t,\\
\bar{X}_0^\e=X_0.
\end{array}
\right.
\de
Since $A_\e$ is Lipschitz continuous, the above equation has a unique strong solution denoted as $\bar{X}_{\cdot}^\e$ (cf. \cite[Theorem 6, P. 249]{p}). Moreover, by the same deduction to that of \cite[Step 1, P. 492]{gz} or \cite[Lemma 6.1]{rwzh}, it holds that 
\be
\lim\limits_{\e\rightarrow0}\mE\left(\sup\limits_{t\in[0,T]}|\bar{X}_t^\e-X_t|^2\right)=0.
\label{barx}
\ee

Besides, note that
\be
\mE\left(\sup\limits_{s\in[0,t]}|X_s^\e-X_s|^2\right)\leq 2\mE\left(\sup\limits_{s\in[0,t]}|X_s^\e-\bar{X}_s^\e|^2\right)+2\mE\left(\sup\limits_{s\in[0,t]}|\bar{X}_s^\e-X_s|^2\right).
\label{xexxebar}
\ee
Thus, we are devoted to estimating $\mE\left(\sup\limits_{s\in[0,t]}|X_s^\e-\bar{X}_s^\e|^2\right)$.

Next, applying the It\^o formula to $|X^\e_t-\bar{X}_t^\e|^2$ for $t\in[0,T]$, by the monotonicity of $A_\e$ we obtain that
\ce
&&|X_t^\e-\bar{X}_t^\e|^2\\
&=&-2\int_0^t\<X_s^\e-\bar{X}_s^\e, A_\e(X_s^\e)-A_\e(\bar{X}_s^\e)\>\dif s+2\int_0^t\<X_s^\e-\bar{X}_s^\e, b(X^\e_s,\sL_{X^\e_s})-b(X_s,\sL_{X_s})\>\dif s\\
&&+2\int_0^t\<X_s^\e-\bar{X}_s^\e, \(\s(X^\e_s,\sL_{X^\e_s})-\s(X_s,\sL_{X_s})\)\dif W_s\>+\int_0^t\|\s(X^\e_s,\sL_{X^\e_s})-\s(X_s,\sL_{X_s})\|^2\dif s\\
&\leq&2\int_0^t\<X_s^\e-X_s, b(X^\e_s,\sL_{X^\e_s})-b(X_s,\sL_{X_s})\>\dif s+2\int_0^t\<X_s-\bar{X}_s^\e, b(X^\e_s,\sL_{X^\e_s})-b(X_s,\sL_{X_s})\>\dif s\\
&&+2\int_0^t\<X_s^\e-\bar{X}_s^\e, \(\s(X^\e_s,\sL_{X^\e_s})-\s(X_s,\sL_{X_s})\)\dif W_s\>+\int_0^t\|\s(X^\e_s,\sL_{X^\e_s})-\s(X_s,\sL_{X_s})\|^2\dif s\\
&=:&I_1(t)+I_2(t)+I_3(t)+I_4(t).
\de
For $I_1$, by (${\bf H}^2_{b}$) it holds that
\ce
\mE\left(\sup\limits_{s\in[0,t]}I_1(s)\right)\leq L_1\mE\int_0^t\left(|X_s^\e-X_s|^2+\mW_2^2(\sL_{X^\e_s},\sL_{X_s})\right)\dif s\leq 2L_1\int_0^t\mE\left(\sup\limits_{r\in[0,s]}|X_r^\e-X_r|^2\right)\dif s.
\de
And we deal with $I_2$. (${\bf H}^1_{b,\s}$) and Lemma \ref{xkes}, \ref{Yeqsomo} imply that
\ce
\mE\left(\sup\limits_{s\in[0,t]}I_2(s)\right)&\leq& 2L_{b,\s}\mE\int_0^t|X_s-\bar{X}_s^\e|(2+|X^\e_s|+\|\sL_{X^\e_s}\|_2+|X_s|+\|\sL_{X_s}\|_2)\dif s\\
&\leq& 2L_{b,\s}\mE\left(\sup\limits_{s\in[0,t]}|X_s-\bar{X}_s^\e|\right)\left(\int_0^t(2+|X^\e_s|+\|\sL_{X^\e_s}\|_2+|X_s|+\|\sL_{X_s}\|_2)\dif s\right)\\
&\leq&C\left(\mE\left(\sup\limits_{s\in[0,t]}|X_s-\bar{X}_s^\e|^2\right)\right)^{1/2}\left(\int_0^t(2+\mE|X^\e_s|^2+\mE|X_s|^2)\dif s\right)^{1/2}\\
&\leq&C\left(\mE\left(\sup\limits_{s\in[0,T]}|X_s-\bar{X}_s^\e|^2\right)\right)^{1/2}.
\de
For $I_3$, from (${\bf H}^2_{\s}$), the BDG inequality and the H\"older inequality, it follows that 
\ce
\mE\left(\sup\limits_{s\in[0,t]}I_3(s)\right)&\leq&C\mE\left(\int_0^t|X_s^\e-\bar{X}_s^\e|^2\|\s(X^\e_s,\sL_{X^\e_s})-\s(X_s,\sL_{X_s})\|^2\dif s\right)^{1/2}\\
&\leq&\frac{1}{2}\mE\left(\sup\limits_{s\in[0,t]}|X_s^\e-\bar{X}_s^\e|^2\right)+C\mE\int_0^t\(|X^\e_s-X_s|^2+\mW_2^2(\sL_{X^\e_s},\sL_{X_s})\)\dif s\\
&\leq&\frac{1}{2}\mE\left(\sup\limits_{s\in[0,t]}|X_s^\e-\bar{X}_s^\e|^2\right)+C\int_0^t\mE\left(\sup\limits_{r\in[0,s]}|X_r^\e-X_r|^2\right)\dif s.
\de
By (${\bf H}^2_{\s}$), we treat $I_4$ to obtain that
\ce
\mE\left(\sup\limits_{s\in[0,t]}I_4(s)\right)\leq L_2\mE\int_0^t\(|X^\e_s-X_s|^2+\mW_2^2(\sL_{X^\e_s},\sL_{X_s})\)\dif s\leq 2L_2\int_0^t\mE\left(\sup\limits_{r\in[0,s]}|X_r^\e-X_r|^2\right)\dif s.
\de
Combining all the above estimates, we have that
\be
\mE\left(\sup\limits_{s\in[0,t]}|X_s^\e-\bar{X}_s^\e|^2\right)\leq C\left(\mE\left(\sup\limits_{s\in[0,T]}|X_s-\bar{X}_s^\e|^2\right)\right)^{1/2}+C\int_0^t\mE\left(\sup\limits_{r\in[0,s]}|X_r^\e-X_r|^2\right)\dif s.
\label{xebarx}
\ee

Finally, (\ref{xexxebar}) (\ref{xebarx}) yield that
\ce
\mE\left(\sup\limits_{s\in[0,t]}|X_s^\e-X_s|^2\right)&\leq&C\left(\mE\left(\sup\limits_{s\in[0,T]}|X_s-\bar{X}_s^\e|^2\right)\right)^{1/2}+2\mE\left(\sup\limits_{s\in[0,T]}|\bar{X}_s^\e-X_s|^2\right)\\
&&+C\int_0^t\mE\left(\sup\limits_{r\in[0,s]}|X_r^\e-X_r|^2\right)\dif s.
\de
By the Gronwall inequality and (\ref{barx}), we obtain the required first result. 

The same deduction to that for the first limit gives the second limit. The proof is complete.
\end{proof}

\bp\label{xnexees}
Assume that $b, \s, b^n, \s^n$ uniformly satisfy (${\bf H}^1_{b,\s}$) (${\bf H}^2_{b}$) (${\bf H}^2_{\s}$) and $X_0=X^n_0, \mE|X_0|^{2p}<\infty$ for any $p>2$. Moreover, (${\bf H}^1_{A^n,A}$) (${\bf H}^2_{A^n,A}$) (${\bf H}_{b^n,b,\s^n,\s}$) hold. Then for any $\e>0$
\ce
\lim\limits_{n\rightarrow\infty}\sup\limits_{t\in[0,T]}\mE|X^{n,\e}_t-X^{\e}_t|^2=0.
\de
\ep
\begin{proof}
First of all, note that
\ce
&&X^{n,\e}_t=X_0-\int_0^tA^n_{\e}(X^{n,\e}_s)\dif s+\int_0^tb^n(X^{n,\e}_s,\sL_{X^{n,\e}_s})\dif s+\int_0^t\sigma^n(X^{n,\e}_s,\sL_{X^{n,\e}_s})\dif W_s,\\
&&X^{\e}_t=X_0-\int_0^tA_{\e}(X^{\e}_s)\dif s+\int_0^tb(X^{\e}_s,\sL_{X^{\e}_s})\dif s+\int_0^t\sigma(X^{\e}_s,\sL_{X^{\e}_s})\dif W_s.
\de
Thus, applying the It\^o formula to $|X^{n,\e}_t-X^{\e}_t|^2$ and taking the expectation on two sides, we obtain that
\be
\mE|X^{n,\e}_t-X^{\e}_t|^2&=&-2\mE\int_0^t\<X^{n,\e}_s-X^{\e}_s,A^n_{\e}(X^{n,\e}_s)-A_{\e}(X^{\e}_s)\>\dif s\no\\
&&+2\mE\int_0^t\<X^{n,\e}_s-X^{\e}_s, b^n(X^{n,\e}_s,\sL_{X^{n,\e}_s})-b(X^{\e}_s,\sL_{X^{\e}_s})\>\dif s\no\\
&&+\mE\int_0^t\|\sigma^n(X^{n,\e}_s,\sL_{X^{n,\e}_s})-\sigma(X^{\e}_s,\sL_{X^{\e}_s})\|^2\dif s\no\\
&=:&J_1(t)+J_2(t)+J_3(t).
\label{j1j2j3}
\ee

For $J_1$, by the monotonicity of $A^n_{\e}$ and the H\"older inequality, it holds that
\be
J_1(t)&=&-2\mE\int_0^t\<X^{n,\e}_s-X^{\e}_s,A^n_{\e}(X^{n,\e}_s)-A^n_{\e}(X^{\e}_s)+A^n_{\e}(X^{\e}_s)-A_{\e}(X^{\e}_s)\>\dif s\no\\
&\leq&-2\mE\int_0^t\<X^{n,\e}_s-X^{\e}_s,A^n_{\e}(X^{\e}_s)-A_{\e}(X^{\e}_s)\>\dif s\no\\
&\leq&\int_0^t\mE|X^{n,\e}_s-X^{\e}_s|^2\dif s+\int_0^t\mE|A^n_{\e}(X^{\e}_s)-A_{\e}(X^{\e}_s)|^2\dif s\no\\
&\leq&\int_0^t\mE|X^{n,\e}_s-X^{\e}_s|^2\dif s+\int_0^t\mE|A^n_{\e}(X^{\e}_s)-A_{\e}(X^{\e}_s)|^2I_{\left\{\sup\limits_{s\in[0,T]}|X^{\e}_s|> R_1\right\}}\dif s\no\\
&&+\int_0^t\mE|A^n_{\e}(X^{\e}_s)-A_{\e}(X^{\e}_s)|^2I_{\left\{\sup\limits_{s\in[0,T]}|X^{\e}_s|\leq R_1\right\}}\dif s,
\label{j1}
\ee
where $R_1>0$ is a constant determined later. And we estimate the second and third terms for the right side of the above inequality. Note that for any $x\in\cD(A)$
\ce
|A^n_{\e}(x)|+|A_{\e}(x)|\leq |(A^n)^\circ(0)|+|A^\circ(0)|+\frac{2}{\e}|x|\leq \left(\g+|A^\circ(0)|+\frac{2}{\e}\right)(1+|x|),
\de
where the Lipschitz continuity of $A^n_{\e}, A_{\e}$ and $|A^n_{\e}(0)|\leq |(A^n)^\circ(0)|\leq \g, |A_{\e}(0)|\leq |A^\circ(0)|$ are used. Thus, the Chebyshev inequality implies that
\ce
&&\int_0^t\mE|A^n_{\e}(X^{\e}_s)-A_{\e}(X^{\e}_s)|^2I_{\left\{\sup\limits_{s\in[0,T]}|X^{\e}_s|> R_1\right\}}\dif s\\
&\leq&2\left(\g+|A^\circ(0)|+\frac{2}{\e}\right)^2\int_0^t\mE(1+|X^{\e}_s|^2)I_{\left\{\sup\limits_{s\in[0,T]}|X^{\e}_s|> R_1\right\}}\dif s\\
&\leq&2\left(\g+|A^\circ(0)|+\frac{2}{\e}\right)^2T\frac{\mE\left(\sup\limits_{s\in[0,T]}|X^{\e}_s|^2+\sup\limits_{s\in[0,T]}|X^{\e}_s|^4\right)}{R_1^2}.
\de
And for any $\d>0$, we take $R_1$ large enough such that 
\be
\int_0^T\mE|A^n_{\e}(X^{\e}_s)-A_{\e}(X^{\e}_s)|^2I_{\{\sup\limits_{s\in[0,T]}|X^{\e}_s|> R_1\}}\dif s\leq \d.
\label{j11}
\ee
Besides, by (${\bf H}^2_{A^n,A}$), it holds that
\ce
\lim\limits_{n\rightarrow\infty}|A^n_{\e}(X^{\e}_s)-A_{\e}(X^{\e}_s)|^2I_{\{\sup\limits_{s\in[0,T]}|X^{\e}_s|\leq R_1\}}=0.
\de
So, by the dominated convergence theorem, we have that
\be
\lim\limits_{n\rightarrow\infty}\int_0^T\mE|A^n_{\e}(X^{\e}_s)-A_{\e}(X^{\e}_s)|^2I_{\{\sup\limits_{s\in[0,T]}|X^{\e}_s|\leq R_1\}}\dif s=0.
\label{j12}
\ee

Next, we deal with $J_2(t)+J_3(t)$. (${\bf H}^2_{b}$) (${\bf H}^2_{\s}$) and the H\"older inequality imply that
\be
&&J_2(t)+J_3(t)\no\\
&=&2\mE\int_0^t\<X^{n,\e}_s-X^{\e}_s, b^n(X^{n,\e}_s,\sL_{X^{n,\e}_s})-b^n(X^{\e}_s,\sL_{X^{\e}_s})+b^n(X^{\e}_s,\sL_{X^{\e}_s})-b(X^{\e}_s,\sL_{X^{\e}_s})\>\dif s\no\\
&&+\mE\int_0^t\|\sigma^n(X^{n,\e}_s,\sL_{X^{n,\e}_s})-\sigma^n(X^{\e}_s,\sL_{X^{\e}_s})+\sigma^n(X^{\e}_s,\sL_{X^{\e}_s})-\sigma(X^{\e}_s,\sL_{X^{\e}_s})\|^2\dif s\no\\
&\leq&(L_1+2L_2)\mE\int_0^t(|X^{n,\e}_s-X^{\e}_s|^2+\mW_2^2(\sL_{X^{n,\e}_s},\sL_{X^{\e}_s}))\dif s+\mE\int_0^t|X^{n,\e}_s-X^{\e}_s|^2\dif s\no\\
&&+\mE\int_0^t|b^n(X^{\e}_s,\sL_{X^{\e}_s})-b(X^{\e}_s,\sL_{X^{\e}_s})|^2\dif s+2\mE\int_0^t\|\sigma^n(X^{\e}_s,\sL_{X^{\e}_s})-\sigma(X^{\e}_s,\sL_{X^{\e}_s})\|^2\dif s\no\\
&\leq&(2L_1+4L_2+1)\mE\int_0^t|X^{n,\e}_s-X^{\e}_s|^2\dif s\no\\
&&+\mE\int_0^t|b^n(X^{\e}_s,\sL_{X^{\e}_s})-b(X^{\e}_s,\sL_{X^{\e}_s})|^2I_{\{\sup\limits_{s\in[0,T]}|X^{\e}_s|> R_2\}}\dif s\no\\
&&+2\mE\int_0^t\|\s^n(X^{\e}_s,\sL_{X^{\e}_s})-\s(X^{\e}_s,\sL_{X^{\e}_s})\|^2I_{\{\sup\limits_{s\in[0,T]}|X^{\e}_s|> R_2\}}\dif s\no\\
&&+\mE\int_0^t|b^n(X^{\e}_s,\sL_{X^{\e}_s})-b(X^{\e}_s,\sL_{X^{\e}_s})|^2I_{\{\sup\limits_{s\in[0,T]}|X^{\e}_s|\leq R_2\}}\dif s\no\\
&&+2\mE\int_0^t\|\s^n(X^{\e}_s,\sL_{X^{\e}_s})-\s(X^{\e}_s,\sL_{X^{\e}_s})\|^2I_{\{\sup\limits_{s\in[0,T]}|X^{\e}_s|\leq R_2\}}\dif s,
\label{j2}
\ee
where $R_2>0$ is a constant determined later. By (${\bf H}^1_{b,\s}$) (${\bf H}_{b^n,b,\s^n,\s}$) and the same deduction to that for (\ref{j11}) (\ref{j12}), we obtain that for $R_2$ large enough
\be
&&\bigg(\mE\int_0^T|b^n(X^{\e}_s,\sL_{X^{\e}_s})-b(X^{\e}_s,\sL_{X^{\e}_s})|^2I_{\{\sup\limits_{s\in[0,T]}|X^{\e}_s|> R_2\}}\dif s\no\\
&&\quad +2\mE\int_0^T\|\s^n(X^{\e}_s,\sL_{X^{\e}_s})-\s(X^{\e}_s,\sL_{X^{\e}_s})\|^2I_{\{\sup\limits_{s\in[0,T]}|X^{\e}_s|> R_2\}}\dif s\bigg)\leq \d,\label{j21}\\
&&\lim\limits_{n\rightarrow\infty}\bigg(\int_0^T\mE|b^n(X^{\e}_s,\sL_{X^{\e}_s})-b(X^{\e}_s,\sL_{X^{\e}_s})|^2I_{\{\sup\limits_{s\in[0,T]}|X^{\e}_s|\leq R_2\}}\dif s\no\\
&&\quad +2\mE\int_0^T\|\s^n(X^{\e}_s,\sL_{X^{\e}_s})-\s(X^{\e}_s,\sL_{X^{\e}_s})\|^2I_{\{\sup\limits_{s\in[0,T]}|X^{\e}_s|\leq R_2\}}\dif s\bigg)=0.
\label{j22}
\ee

Finally, combining (\ref{j1}) (\ref{j11}) (\ref{j2}) (\ref{j21}) with (\ref{j1j2j3}), we have that
\ce
&&\sup\limits_{s\in[0,t]}\mE|X^{n,\e}_s-X^{\e}_s|^2\\
&\leq& C\int_0^t\sup\limits_{r\in[0,s]}\mE|X^{n,\e}_r-X^{\e}_r|^2\dif s+2\d+\int_0^T\mE|A^n_{\e}(X^{\e}_s)-A_{\e}(X^{\e}_s)|^2I_{\{\sup\limits_{s\in[0,T]}|X^{\e}_s|\leq R_1\}}\dif s\\
&&+\int_0^T\mE|b^n(X^{\e}_s,\sL_{X^{\e}_s})-b(X^{\e}_s,\sL_{X^{\e}_s})|^2I_{\{\sup\limits_{s\in[0,T]}|X^{\e}_s|\leq R_2\}}\dif s\\
&&+2\mE\int_0^T|\s^n(X^{\e}_s,\sL_{X^{\e}_s})-\s(X^{\e}_s,\sL_{X^{\e}_s})|^2I_{\{\sup\limits_{s\in[0,T]}|X^{\e}_s|\leq R_2\}}\dif s.
\de
So, the Gronwall inequality and (\ref{j12}) (\ref{j22}) yields the required result. The proof is complete.
\end{proof}

Now, it is the position to prove Theorem \ref{constrsol}.

{\bf Proof of Theorem \ref{constrsol}.} First of all, note that for any $\e>0$
\ce
&&\sup\limits_{t\in[0,T]}\mE|X^{n}_t-X_t|^2\\
&\leq& \sup\limits_{t\in[0,T]}\mE|X^{n}_t-X^{n,\e}_t|^2+\sup\limits_{t\in[0,T]}\mE|X^{n,\e}_t-X^{\e}_t|^2+\sup\limits_{t\in[0,T]}\mE|X^{\e}_t-X_t|^2\\
&\leq&\mE\left(\sup\limits_{t\in[0,T]}|X^{n}_t-X^{n,\e}_t|^2\right)+\sup\limits_{t\in[0,T]}\mE|X^{n,\e}_t-X^{\e}_t|^2+\mE\left(\sup\limits_{t\in[0,T]}|X^{\e}_t-X_t|^2\right).
\de
Thus, as $n\rightarrow\infty$, by Proposition \ref{xnexees} it holds that
\ce
\lim\limits_{n\rightarrow\infty}\sup\limits_{t\in[0,T]}\mE|X^{n}_t-X_t|^2\leq \sup\limits_{n}\mE\left(\sup\limits_{t\in[0,T]}|X^{n}_t-X^{n,\e}_t|^2\right)+\mE\left(\sup\limits_{t\in[0,T]}|X^{\e}_t-X_t|^2\right),
\de
which together with Proposition \ref{xexcon} completes the proof.

\section{Convergence of invariant measures}\label{coninvmease}

In the section, we investigate the convergence of invariant measures for a sequence of multivalued McKean-Vlasov SDEs. 

Consider Eq.(\ref{eq1}) and Eq.(\ref{eqn}), i.e.
\ce
&&\dif X_t\in \ -A(X_t)\dif t+b(X_t,\sL_{X_t})\dif t+\sigma(X_t,\sL_{X_t})\dif W_t,\\
&&\dif X^n_t\in \ -A^n(X^n_t)\dif t+b^n(X^n_t,\sL_{X^n_t})\dif t+\sigma^n(X^n_t,\sL_{X^n_t})\dif W_t.
\de
When $A^n, A$ satisfy (${\bf H}_{A}$), and $b, \s, b^n, \s^n$ uniformly satisfy (${\bf H}^1_{b,\s}$) (${\bf H}^2_{\s}$) (${\bf H}^{2'}_{b,\s}$), by Theorem \ref{erg2}
Eq.(\ref{eq1}) with $X_0\in\overline{\cD(A)}, \mE|X_0|^2<\infty$ and Eq.(\ref{eqn}) with $X^n_0\in\overline{\cD(A^n)}, \mE|X^n_0|^2<\infty$ have unique invariant probability measures $\mu_{\infty},\mu^n_{\infty}\in\cP_{2}(\overline{\cD(A)})$, respectively. About $\mu_{\infty},\mu^n_{\infty}$, we have the following conclusion which is the main result in this section.

\bt\label{coninvmea}
Assume that $b, \s, b^n, \s^n$ uniformly satisfy (${\bf H}^1_{b,\s}$) (${\bf H}^2_{\s}$) (${\bf H}^{2'}_{b,\s}$) and $X_0=X^n_0, \mE|X_0|^{2p}<\infty$ for any $p>2$. Moreover, (${\bf H}^1_{A^n,A}$) (${\bf H}^2_{A^n,A}$) (${\bf H}_{b^n,b,\s^n,\s}$) hold. Then as $n\rightarrow\infty$
\ce
\mu^n_{\infty}\overset{w}{\longrightarrow}\mu_{\infty}.
\de
\et
\begin{proof}
{\bf Step 1.} We prove that $\{\mu^n_{\infty}, n\in\mN\}$ is tight.

First of all, we consider Eq.(\ref{eqn}) with $\sL_{X_0^n}=\mu^n_{\infty}$ and denote its solution by $(X_{\cdot}^{n,\infty}, K_{\cdot}^{n,\infty})$. Then by applying the It\^o formula to $|X^{n,\infty}_t|^2$, it holds that for any $y\in A^n(0)$
\ce
|X^{n,\infty}_t|^2&=&|X^{n,\infty}_0|^2-\int_0^t2\<X^{n,\infty}_s,\dif K^{n,\infty}_s\>+\int_0^t2\<X^{n,\infty}_s,b^{n}(X^{n,\infty}_s,\sL_{X^{n,\infty}_s})\>\dif s\\
&&+\int_0^t2\<X^{n,\infty}_s,\s^{n}(X^{n,\infty}_s,\sL_{X^{n,\infty}_s})\dif W_s\>+\int_0^t\|\s^{n}(X^{n,\infty}_s,\sL_{X^{n,\infty}_s})\|^2\dif s\\
&\leq&|X^{n,\infty}_0|^2+2|y|\int_0^t|X^{n,\infty}_s|\dif s+\int_0^t2\<X^{n,\infty}_s,b^{n}(X^{n,\infty}_s,\sL_{X^{n,\infty}_s})\>\dif s\\
&&+\int_0^t2\<X^{n,\infty}_s,\s^{n}(X^{n,\infty}_s,\sL_{X^{n,\infty}_s})\dif W_s\>+\int_0^t\|\s^{n}(X^{n,\infty}_s,\sL_{X^{n,\infty}_s})\|^2\dif s,
\de
where Lemma \ref{equi} is used. And taking the expectation on two sides of the above inequality, by (\ref{ergcon}) we get that
\ce
\mE|X^{n,\infty}_t|^2&\leq& \mE|X^{n,\infty}_0|^2+2|y|\int_0^t\mE|X^{n,\infty}_s|\dif s\\
&&+\int_0^t\(C+(L_3+L_2)\|\sL_{X^{n,\infty}_s}\|^2_2-(L_4-L_2-\frac{\lambda-2L_2}{4})\mE|X^{n,\infty}_s|^2\)\dif s.
\de
Note that $P_t^{n*}\mu^n_{\infty}=\mu^n_{\infty}$, where $P_t^{n*}\mu^n_{\infty}:=\sL_{X^{n,\infty}_t}$ for $\sL_{X_0^{n,\infty}}=\mu^n_{\infty}$. Thus, it holds that
\ce
\frac{3(\lambda-2L_2)}{4}\int_{\overline{\cD(A)}}|x|^2\mu^n_{\infty}(\dif x)\leq 2|y|\int_{\overline{\cD(A)}}|x|\mu^n_{\infty}(\dif x)+C.
\de
Besides, for any $r>\frac{4(2|y|+1)}{3(\lambda-2L_2)}$, we know that 
$$
\frac{|x|}{r}\leq 1+\frac{|x|^2}{r^2}, \quad \forall x\in\mR^d,
$$
which yields that
\ce
\frac{3(\lambda-2L_2)}{4}\int_{\overline{\cD(A)}}|x|^2\mu^n_{\infty}(\dif x)&\leq& 2|y|\int_{\overline{\cD(A)}}\(r+\frac{|x|^2}{r}\)\mu^n_{\infty}(\dif x)+C\\
&=&2|y|r+C+\frac{2|y|}{r}\int_{\overline{\cD(A)}}|x|^2\mu^n_{\infty}(\dif x),
\de
and furthermore
\be
\sup\limits_{n}\int_{\overline{\cD(A)}}|x|^2\mu^n_{\infty}(\dif x)\leq C,
\label{unibou}
\ee
where $C>0$ is independent of $n$. Set $B_R:=\{x\in\mR^d: |x|\leq R\}$, and then
\ce
\sup\limits_{n}\mu^n_{\infty}(B^c_R\cap\overline{\cD(A)})\leq \frac{1}{R^2}\int_{B^c_R\cap\overline{\cD(A)}}|x|^2\mu^n_{\infty}(\dif x)\leq \frac{C}{R^2},
\de
which implies that $\{\mu^n_{\infty}, n\in\mN\}$ is tight.

{\bf Step 2.} We prove that the limit of $\{\mu^n_{\infty}, n\in\mN\}$ is $\mu_{\infty}$.

Since $\{\mu^n_{\infty}, n\in\mN\}$ is tight, there exist a subsequence still denoted as $\{\mu^n_{\infty}, n\in\mN\}$ and a probability measure $\bar{\mu}_{\infty}\in\cP_{2}(\overline{\cD(A)})$ such that as $n\rightarrow\infty$
\be
\mu^n_{\infty}\overset{w}{\longrightarrow}\bar{\mu}_{\infty}.
\label{munweakb}
\ee
To prove that $\bar{\mu}_{\infty}=\mu_{\infty}$, noticing that $\mu_{\infty}$ is the unique invariant probability measure of $\{P_{t}^{*}, t\geq 0\}$, we only need to show that $\bar{\mu}_{\infty}$ is an invariant probability measure of $\{P_{t}^{*}, t\geq 0\}$, i.e. $P_{t}^{*}\bar{\mu}_{\infty}=\bar{\mu}_{\infty}$ for all $t\geq 0$.

Next, for any $\varphi\in C_{lip,b}(\overline{\cD(A)})$ with the Lipschitz constant $L_\varphi$, it holds that
\ce
\left|\int_{\overline{\cD(A)}}\varphi\dif \bar{\mu}_{\infty}-\int_{\overline{\cD(A)}}\varphi\dif P_{t}^{*}\bar{\mu}_{\infty}\right|&\leq& \left|\int_{\overline{\cD(A)}}\varphi\dif \bar{\mu}_{\infty}-\int_{\overline{\cD(A)}}\varphi\dif P_{t}^{n*}\mu^n_{\infty}\right|\\
&&+\left|\int_{\overline{\cD(A)}}\varphi\dif P_{t}^{n*}\mu^n_{\infty}-\int_{\overline{\cD(A)}}\varphi\dif P_{t}^{n*}\bar{\mu}_{\infty}\right|\\
&&+\left|\int_{\overline{\cD(A)}}\varphi\dif P_{t}^{n*}\bar{\mu}_{\infty}-\int_{\overline{\cD(A)}}\varphi\dif P_{t}^{*}\bar{\mu}_{\infty}\right|\\
&=:&I_1+I_2+I_3.
\de
For $I_1$, noticing $P_t^{n*}\mu^n_{\infty}=\mu^n_{\infty}$, by (\ref{munweakb}) we have that
\be
\lim\limits_{n\rightarrow\infty}I_1=0.
\label{i1}
\ee
We deal with $I_2$. On one side, by the same deduction to that for (\ref{xynumu}), it holds that
\be
I_2&=&\left|\mE \varphi(X_t^{n,\mu^n_{\infty}})-\mE \varphi(X_t^{n,\bar{\mu}_{\infty}})\right|\leq \mE|\varphi(X_t^{n,\mu^n_{\infty}})-\varphi(X_t^{n,\bar{\mu}_{\infty}})|\leq L_\varphi\mE|X_t^{n,\mu^n_{\infty}}-X_t^{n,\bar{\mu}_{\infty}}|\no\\
&\leq&L_\varphi(\mE|X_t^{n,\mu^n_{\infty}}-X_t^{n,\bar{\mu}_{\infty}}|^2)^{1/2}\leq L_\varphi\mW_2(\mu^n_{\infty},\bar{\mu}_{\infty}).
\label{i2w2}
\ee
On the other side, since $\mu^n_{\infty}\overset{w}{\longrightarrow}\bar{\mu}_{\infty}$ as $n\rightarrow\infty$, we know that there exist a probability space $(\tilde{\Omega},\tilde{\mathscr{F}},\tilde{\mP})$ and $\overline{\cD(A)}$-valued random variables $\xi^n$, $\xi$ on it satisfying $\xi^n\stackrel{a.s.}\longrightarrow\xi$ as $n\rightarrow\infty$. And the estimate (\ref{unibou}) assures that $\sup\limits_{n}\tilde{\mE}|\xi^n|^2\leq C$, where $\tilde{\mE}$ stands for the expectation under the probability measure $\tilde{\mP}$. So, the dominated convergence theorem implies that
$$
\lim\limits_{n\rightarrow\infty}\tilde{\mE}|\xi^n-\xi|^2=0,
$$
which together with $\mW_2(\mu^n_{\infty},\bar{\mu}_{\infty})\leq(\tilde{\mE}|\xi^n-\xi|^2)^{1/2}$ yields that
\be
\lim\limits_{n\rightarrow\infty}\mW_2(\mu^n_{\infty},\bar{\mu}_{\infty})=0.
\label{w2ze}
\ee
Thus, (\ref{i2w2}) and (\ref{w2ze}) implies that
\be
\lim\limits_{n\rightarrow\infty}I_2=0.
\label{i2}
\ee
To $I_3$, by simple calculation, one can obtain that
\ce
I_3&=&\left|\mE \varphi(X_t^{n,\bar{\mu}_{\infty}})-\mE \varphi(X_t^{\bar{\mu}_{\infty}})\right|\leq\mE\left|\varphi(X_t^{n,\bar{\mu}_{\infty}})-\varphi(X_t^{\bar{\mu}_{\infty}})\right|\leq L_\varphi\mE\left|X_t^{n,\bar{\mu}_{\infty}}-X_t^{\bar{\mu}_{\infty}}\right|\\
&\leq& L_\varphi\(\mE|X_t^{n,\bar{\mu}_{\infty}}-X_t^{\bar{\mu}_{\infty}}|^2\)^{1/2},
\de
which together with Theorem \ref{constrsol} yields that
\be
\lim\limits_{n\rightarrow\infty}I_3=0.
\label{i3}
\ee

Finally, combining (\ref{i1}) (\ref{i2}) (\ref{i3}), we conclude that 
$$
\rho(\bar{\mu}_{\infty},P_{t}^{*}\bar{\mu}_{\infty})=\sup\limits_{\parallel{\varphi}\parallel_{C_{lip,b}(\overline{\cD(A)})}\leq1}\left|\int_{\overline{\cD(A)}}\varphi\dif \bar{\mu}_{\infty}-\int_{\overline{\cD(A)}}\varphi\dif P_{t}^{*}\bar{\mu}_{\infty}\right|=0,
$$
which implies that $P_{t}^{*}\bar{\mu}_{\infty}=\bar{\mu}_{\infty}$ for all $t\geq 0$. The proof is complete.
\end{proof}

\end{document}